# Discrete Optimal Designs for Distributed Energy Systems with Nonconvex Multiphase Optimal Power Flow

*Ishanki De Mel[a], Oleksiy V. Klymenko[a], Michael Short[a*]*

[a]Department of Chemical and Process Engineering, University of Surrey, Guildford, United Kingdom

GU2 7XH

*m.short@surrey.ac.uk

## Abstract

The optimal selection, sizing, and location of small-scale technologies within a grid-connected distributed energy system (DES) can contribute to reducing carbon emissions, consumer costs, and network imbalances. There is a significant lack of studies on how DES designs, especially those with electrified heating systems, impact unbalanced low-voltage distribution networks to which most DES are connected. This is the first study to present an optimisation framework for obtaining discrete technology sizing and selection for grid-connected DES design, while simultaneously considering multiphase optimal power flow (MOPF) constraints to accurately represent unbalanced low-voltage distribution networks. An algorithm is developed to solve the resulting Mixed-Integer Nonlinear Programming (MINLP) formulation. It employs a decomposition based on Mixed-Integer Linear Programming (MILP) and Nonlinear Programming (NLP), and utilises integer cuts and complementarity reformulations to obtain discrete designs that are also feasible with respect to the network constraints. A heuristic modification to the original algorithm is also proposed to improve computational speed. Improved formulations for selecting feasible combinations of air source heat pumps (ASHPs) and hot water storage tanks are also presented. Two networks of varying size are used to test the optimisation methods. Designs with electrified heating (ASHPs and tanks) are compared to those with conventional gas boilers. The algorithms outperform the existing state-of-the-art commercial MINLP solver, which fails to find any solutions in two instances. While feasible solutions were obtained for all cases, convergence was not achieved for all, especially for those involving the larger network. Where converged, the algorithm with the heuristic modification has achieved results up to 70% faster than the original algorithm. Results for case studies suggest that including ASHPs can support up to 16% higher renewable generation capacity compared to gas boilers, albeit with higher ASHP investment costs, as local generation and consumption minimises network violations associated with excess power export. The optimisation framework and results can be used to inform stakeholders such as policy-makers and network operators, to increase renewable energy capacity and aid the decarbonisation of domestic heating systems.

**Keywords:** Distributed Energy, Technology Selection, Discrete Sizing, Nonlinear, Multiphase

## Abbreviations

| | |
|---|---|
| AC | Alternating Current |
| ASHP | Air Source Heat Pumps |
| COP | Coefficient of Performance |
| CR | Complementarity reformulations |
| CR-H | Complementarity reformulations with a heuristic modification |
| DER | Distributed Energy Resource |
| DES | Distributed Energy System |
| HW | Hot Water |
| LV | Low Voltage |
| MILP | Mixed-Integer Linear Programming |
| MINLP | Mixed-Integer Nonlinear Programming |
| MOPF | Multiphase Optimal Power Flow |
| NLP | Nonlinear Programming |
| OPF | Optimal Power Flow (typically balanced power flow formulations) |
| PA | Proposed overall algorithm with Algorithm $CR$ embedded. |
| PA-H | Proposed algorithm with Algorithm $CR$-$H$ embedded (includes a heuristic modification) |
| PV | Photovoltaics |

SBB    Simple Branch and Bound (a commercial MINLP solver)

SDP    Semidefinite Programming

## 1. Introduction

Distributed Energy Systems (DES) that operate small-scale renewable energy resources for electricity generation have been increasingly integrated into existing electric power networks [1]. This is primarily due to the efforts to reduce anthropogenic emissions related to energy generation and consumption [2]. DES also typically consist of other distributed energy resources (DERs), such as for energy storage and heat generation, that are located near the consumers [3]. The increasing proportion of electricity from renewable energy resources has also contributed to a focus on electricity-consuming heat generation technologies, such as air source heat pumps (ASHPs), as low-carbon alternatives to conventional fossil fuel-based technologies [4]. DES with renewable-based electricity and heat generation can therefore play a pivotal role in decarbonisation. The "design" of a DES, which involves determining capacities, locations, and types of DERs to be included within a DES [5], has been widely modelled and solved as a mathematical optimisation problem. The design decisions are often best represented using discrete variables, as these enable the optimal selection of technologies from a variety of commercially available units and capacities. DES design models also consider operational constraints and continuous variables to ensure that the designs proposed can meet the consumers' average demands, while considering seasonal and meteorological variations.

Many notable studies on DES design and operation have excluded nonlinear and nonconvex constraints related to alternating current (AC) power networks [6], to which DES are connected. These exclusions have significantly simplified the design problem and enabled the use of Mixed-Integer Linear Programming (MILP) techniques to obtain fast and globally optimal solutions. However, excluding such constraints could produce ill-suited designs for grid-connected DES and contribute to higher operational costs, especially as recent research has shown that DERs can impact grid stability and resilience [7]. More detailed models encapsulating these nonconvex optimal power flow (OPF) constraints have been recently developed, which have been labelled as DES-OPF models in a previous literature review by the same authors [8]. The inclusion of nonconvex power flow constraints within DES design optimisation results in a large-scale Mixed-Integer Nonlinear Programming (MINLP) model that is extremely challenging to solve. The seminal studies in this area, such as Morvaj et al.

[6], Mashayekh et al. [9], and Sfikas et al. [10] have used a variety of techniques and methods to simplify and solve these large-scale models, including linear approximations and finely-tuned metaheuristics. These studies have also assumed that the power flows in these networks are balanced. While this assumption is valid for transmission networks that operate at high voltages over long distances, it does not typically hold for low-voltage distribution networks that have a low reactance-to-resistance ($X/R$) ratio [11]. The presence of single-phase lines connecting consumers to the networks, such as those found in European radial networks, make them inherently unbalanced [12]. Previous work by the same authors has shown that the balanced power flow assumption on unbalanced three-phase networks can predict significantly lower renewable energy capacities across the network [13].

Multiphase Optimal Power Flow (MOPF) refers to a developing class of models that can consider power flows at each phase in multiphase AC distribution networks, to account for unbalanced conditions at steady-state operation [14]. These models contain nonconvex constraints that are challenging to linearise and/or convexify. Some convex approximations have been proposed, which include fixed-point or Taylor series-based linearisations [15], [16], iterative backward-forward methods [17] including those associated with LinDistFlow [18]–[20], and semidefinite programming (SDP) relaxations [21], [22]. Unfortunately, these methods have several flaws that make them ill-suited for integration within DES design models. Assumptions for complex voltage (or voltage phasors) neglect phase shifts, such as those resulting from step-down transformers. Some of these methods also require prior knowledge of the design to predict feasible initial points using specialist simulation software or nonlinear power flow models, which is not usually available to modellers when designing a DES. Note that DES design models are different from DES operational or control models (for e.g. model predictive control [23]), as the latter models exclude design decisions by using a fixed design to focus the assessment on operational performance, often on shorter temporal scales. Some notable studies in this area include Thomas et al. [24] and Jin et al. [25] where OpenDSS [26] (a simulation tool) is used to assess network unbalances, and Morstyn et al. [27] that uses Z-bus calculations [28] for three-phase power flow analysis.

As a result of these challenges, very few studies have attempted to incorporate nonconvex MOPF in optimisation frameworks for DES design, which include discrete decision-making. Dunham et al. [29] propose a method for determining the design of DES with an iterative and linearised multiphase power flow algorithm proposed by Bernstein and Dall'anese [15], but do not consider how the DES designs could influence network imbalances as balanced loads are assumed. Rehman et al. [30] test a framework for designing a DES using HOMER, without incorporating any power flow constraints, followed by a MATLAB/Simulink analysis of multiphase power flows and control. The influence of the MOPF constraints on the DES design has not been investigated in either of these studies. De Mel et al. [13] successfully integrate MOPF constraints to DES design, but propose a method where key design variables for technologies must remain continuous for the nonlinear constraints to influence the design. Unfortunately, this limits the inclusion of discrete technologies such as ASHPs, for which parametric descriptions cannot be made continuous and applicable to a range of capacities. Many DES studies in the past have included electricity-consuming heating technologies such as ASHPs in design frameworks. However, a framework encompassing optimal design of these technologies when connected to unbalanced networks has not been proposed before. This is an important gap to be addressed, as previous Monte Carlo simulations for a rural Belgian low-voltage network have indicated that a 30% heat pump penetration can cause significant voltage violations [31]. This can impair present efforts to increase heat pump uptake and decarbonise heating systems. Therefore, it is critical to ensure feasibility of the designs with respect to the nonconvex constraints of the underlying network, as opposed to finding globally optimal solutions with approximations that do not accurately represent these networks.

To summarise, the fundamental question of how grid-connected DES, which include electrified heating systems, can be better designed to avoid violating multiphase network constraints and exacerbating network imbalances has not yet been adequately investigated.

To the authors' best knowledge, there is no optimisation framework that captures discrete decision-making for DES design in the presence of nonconvex MOPF constraints. This includes the optimal design of electricity-consuming heat technologies which are connected to existing unbalanced power networks, as opposed to conventional fossil fuel-based technologies that are uncoupled from the electricity network. The study aims to address these gaps by exploring how better designs can be obtained for DES that are connected to unbalanced distribution networks, and the impacts of including electricity-consuming heat technologies on local generation, consumption, and the underlying networks. To enable this investigation, it is essential to develop reliable and appropriate solution methods to solve a large-scale MINLP and obtain discrete design decisions.

The main contribution of this study is a novel optimisation algorithm based on existing mathematical programming (exact) approaches and techniques to address all the gaps identified above. Note that the novelty of this study lies its application, which is obtaining discrete DES design decisions while considering the power flows of multiphase AC distribution networks to which these DES are connected, and the novel use of existing methods to obtain feasible and locally optimal results. A heuristic modification to the proposed algorithm is also presented and tested, to reduce computational expense and aid faster decision-making.

Improved discrete formulations for ASHP selection in combination with heat storage selection are also presented, which is another contribution of this study. Focusing on previous optimisation models (as opposed to simulation-based frameworks [32]), most studies assume continuous sizing for both heat pumps and hot water tanks [33]–[36]. These studies can recommend capacities that do not align with commercially available discrete units. Other studies that consider discrete sizing [37]–[39] assume that each ASHP and hot water tank capacity combination is feasible, which is usually not the case, as made evident by manufacturers' specifications [40]. This study overcomes these limitations by considering combined ASHP-tank selection for commercially available units in the improved formulation, ensuring that the chosen units can feasibly operate together in practice. The presented formulations also consider Coefficient of Performance (COP) and ASHP capacity as functions of ambient temperature (as per the data provided by manufacturers) [41]. Note that past DES design studies have often oversimplified the latter by assuming a fixed COP and a fixed maximum capacity for all ambient temperatures, either for each season [42] or for all seasons [43]–[45].

The proposed solution methods and formulations for discrete technology selection are presented in Section 2. The solution methods are tested and compared in Section 3, alongside efforts to solve the formulation directly using a commercial MINLP solver. The performance of the proposed methods and their limitations are also discussed.

## 2. Methodology

The focus of the methodology is on obtaining locally optimal solutions that are feasible with respect to nonconvex MOPF constraints, as opposed to "globally optimal" solutions of simplified MILP models that potentially violate MOPF constraints.

This study exploits the structure of the DES-MOPF design problem to decompose and solve the overall MINLP using commercial large-scale solvers, as done in previous work by the same authors [13]. To summarise the previously proposed decomposition, the MILP formulation captures all the linear constraints and associated mixed-integer and continuous variables for DES design, thus providing a lower bound for the overall problem. The NLP formulation captures all the linear and nonlinear constraints, including the nonconvex MOPF constraints, and all continuous variables. This allows the MINLP to be solved first as an MILP, followed by an NLP. However, an NLP formulation does not consider any binary decision variables. In DES design optimisation, both design and operational binary variables exist. As proposed in [13], nonlinear complementarity reformulations in the NLP can replace the linear operational constraints with binary variables that are found in the MILP. Unfortunately, discrete DES design decisions cannot be described using complementarity constraints, which make such NLP reformulations unsuitable for this set of binary decisions. In the past, these discrete decisions have been determined by the MILP and subsequently treated as parameters in the NLP. While this is an efficient decomposition, the nonconvex MOPF constraints have limited influence on the DES design, especially if many DERs are represented using discrete or binary variables. The inability to explore different discrete combinations for the DES design, which may be feasible and locally optimal with respect to nonconvex MOPF constraints in the NLP, is a significant limitation.

The solution method proposed in this study addresses this limitation using integer cuts [46] to explore the discrete combinations for DES design while retaining the benefits of using such a decomposition. The proposed algorithm is presented first in Section 2.1. This includes the original complementarity reformulation algorithm proposed in [13], which helps reduce the number of integer cuts required. A heuristic modification of this algorithm is also proposed in this section, to improve computational speed. The mathematical formulation presented in Section 2.2 builds on previous work by the same authors [13]. Only the modifications relevant to discrete decision-making in the DES design are presented in this paper, and the reader is directed to the original study for the complete formulation, including the verified MOPF constraints.

## 2.1 The proposed solution methods

As mentioned previously, the proposed solution method relies on a lower bound $LB$ provided by the MILP, for which the general form is presented in Eq. (1) below. Note that $LB$ represents the objective function value of the MILP. The latter contains the subset of linear constraints $g(x, y_0, z)$ from the overall MINLP optimisation problem, which include Direct Current (DC) linear approximations for power flow. The continuous variables for DES design denoted by $x$, and $y_0$ represent those for DES DC power flow. Binary variables are indicated as $z$.

$$\min \ LB = f(x, y_0, z) \tag{1a}$$
$$s.t. \ g(x, y, z) \leq 0 \tag{1b}$$
$$x \in \mathbb{R}^m, y \in \mathbb{R}^n, z \in \{0,1\}^p \tag{1c}$$

The NLP, which is presented in general form in Eq. (2), contains nonlinear complementarity reformulations and nonconvex MOPF constraints, indicated by $h(x, y, \widehat{z^r})$. Solving the NLP provides a higher objective value to that of the MILP and is therefore the upper bound $UB$. All binary decisions previously obtained from the MILP can be treated as parameters $\hat{z}$ in the NLP. Alternatively, the complementarity reformulations previously proposed by the same authors [13] reduce the number of binary variables treated as parameters, thus enhancing the number of degrees of freedom available to the NLP. The binary variables remaining after these reformulations are denoted as $z^r \subset z \in \{0,1\}^p$, which are treated as fixed parameters $\widehat{z^r}$ in the NLP formulation. The continuous variables $y$ capture those associated with the nonlinear MOPF constraints.

$$\min \ UB = f(x, z) \tag{2a}$$
$$s.t. \ h(x, y, \widehat{z^r}) \leq 0 \tag{2b}$$
$$x \in \mathbb{R}^m, y \in \mathbb{R}^n \tag{2c}$$

As the NLP is incapable of determining binary decisions, an integer cut is introduced in the MILP to remove the current combination of binary solutions and explore different binary topologies with respect to the NLP [46]:

$$\sum_{\kappa \in B_0} z_\kappa^r + \sum_{\kappa \in B_1} z_\kappa^r - 1 \geq 1 \qquad (1d)$$

where $B_0$ represents the subset of unselected binary decisions, while $B_1$ represents the subset of selected binary decisions [46].

When solving the MILP and NLP iteratively, the lower bound from the stable and convex MILP can be expected to monotonically increase (or not decrease, if symmetrical solutions exist) with each integer cut. However, the upper bound from the NLP may fluctuate due to the presence of nonconvex MOPF constraints (while always remaining higher than the current $LB$). This is portrayed in Figure 1. The lowest upper bound produced by the NLP, i.e., the best available solution to the overall MINLP, is denoted by $LUB$. If the upper bound $UB$ at any iteration is less than the incumbent $LUB$, then the new $LUB$ is equal to the current $UB$. As the MILP continues to explore the binary topology with the aid of integer cuts, the gap between $LB$ and $LUB$ can be expected to decrease. If the lower bound $LB$ exceeds the value of the incumbent $LUB$, it is guaranteed to subsequently produce an upper bound that is greater than the existing $LUB$. At this point, no further iterations are required as no better $UB$ can be obtained, and the overall algorithm has therefore converged. Given sufficient time, the MILP can return an infeasible solution if no further binary combinations exist, i.e., all possible binary combinations have already been assessed. If so, $LUB$ can be returned as the final solution, and this too can be considered a convergence criterion.

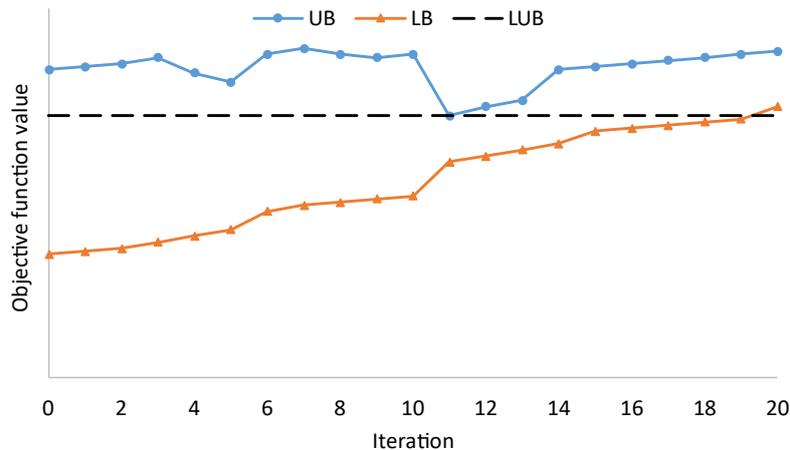

*Figure 1. Example behaviour of the objective functions obtained from the MILP (LB) and NLP (UB) with iterative integer cuts. LUB is the lowest upper bound.*

The overall algorithm encompassing this iterative procedure and convergence criteria is presented in Figure 2. Note that $\lambda$ represents the iteration number or counter. A secondary algorithm for solving nonlinear complementarity reformulations (which eliminates operational binary variables and provides greater degrees of freedom to the NLP) is embedded within this algorithm [13]. This is presented as Algorithm $CR$ in Figure 3. Note that this, too, is an iterative algorithm where $\lambda^{CR}$ represents the iteration number or counter for this algorithm.

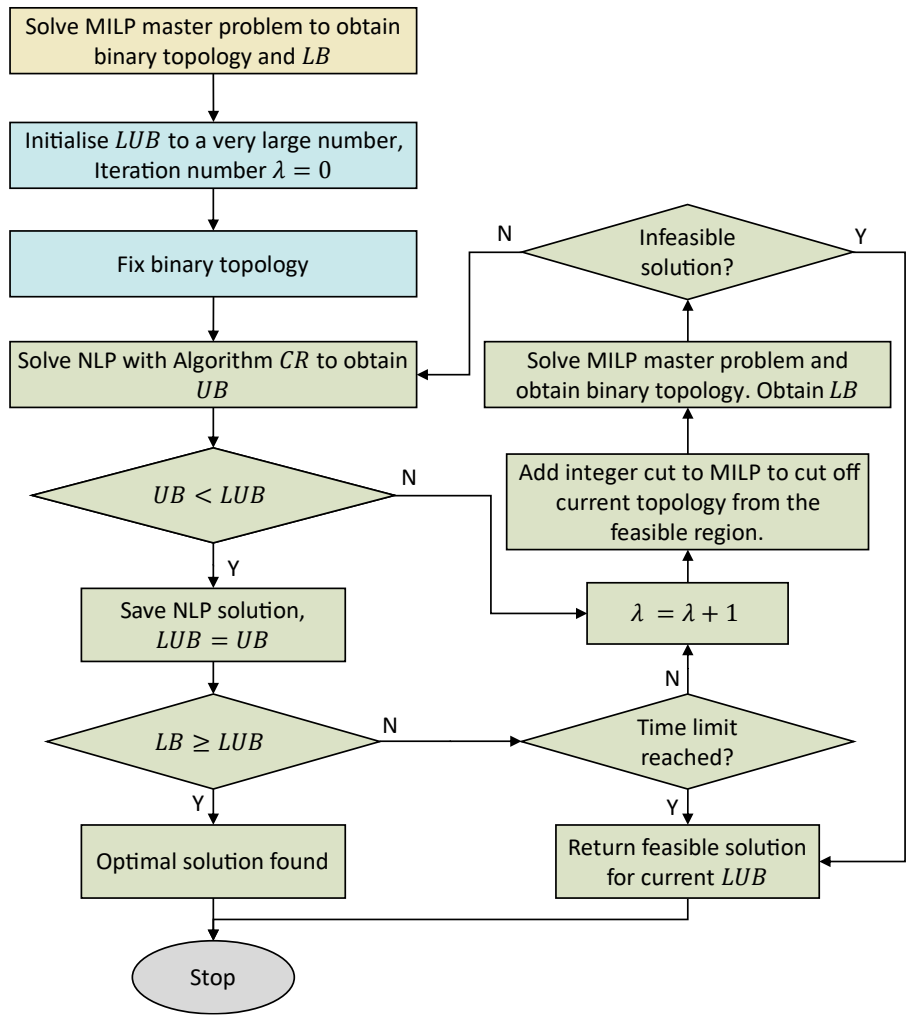

*Figure 2. Proposed overall algorithm with integer cuts.*

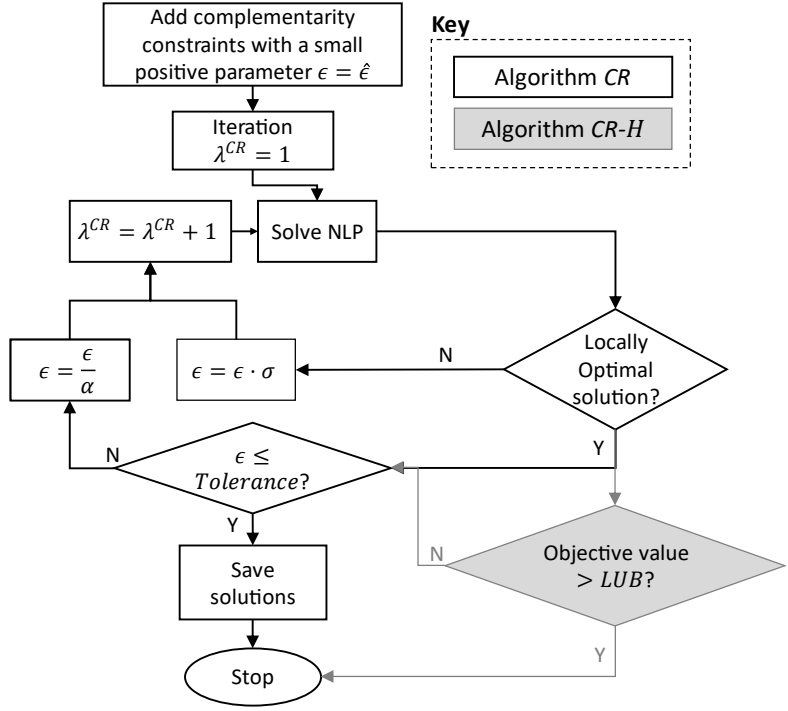

*Figure 3. Algorithms for complementarity reformulations.*

As Algorithm $CR$ is an iterative algorithm for solving the NLP with greater degrees of freedom, computational expense can significantly increase in combination with the overall algorithm. To minimise this computational expense, a heuristic modification to Algorithm $CR$ is proposed, presented as Algorithm $CR\text{-}H$ in Figure 3. This modified algorithm can replace Algorithm $CR$ in the overall algorithm (Figure 2).

The heuristic modification in $CR\text{-}H$ reduces solver time by introducing an early termination condition to the original complementarity reformulations algorithm $CR$. The complementarity reformulations introduce a small positive parameter $\epsilon$ to regularise (relax) the numerically unstable nonlinear formulation shown in Eq. (3) to that in Eq. (4):

$$x \cdot y = 0 \qquad (3)$$
$$x \cdot y \leq \epsilon \qquad (4)$$

This value is iteratively decreased to a value close to zero, such that the complementarity condition is met. In the heuristic modification, following a check in solver status for local optimality, the algorithm is terminated if the objective function value returned by the relaxed NLP is greater than the incumbent lowest upper bound $LUB$. This is because the solver is highly unlikely to find any better objective values (where the objective value is less than the current $LUB$) when the reformulated complementarity constraints are tightened with smaller $\epsilon$ values in subsequent iterations. This condition is bypassed if local optimality is not retained, as this indicates instability whereby $\epsilon$ should be increased to find a locally optimal solution. If the solution returns an objective value less than or equal to $LUB$, this indicates that the MILP has found an integer solution which has improved the NLP upper bound, and the algorithm can proceed to find a solution that meets the complementarity condition.

The heuristic condition imposed in Algorithm $CR\text{-}H$ has its limitations. Considering the nonconvexity of the NLP due to the presence of MOPF constraints, there is a possibility for the NLP to find better or worse locally optimal objective function values in subsequent iterations of the complementarity reformulations algorithm, when compared to the initial prediction at the first iteration. Due to this, the efficacy of the $CR\text{-}H$ algorithm within the overall algorithm is tested by comparing its results to that of the $CR$ algorithm within the overall algorithm.

## 2.2 Mathematical formulation

The objective function minimises total annualised cost $TAC$, which includes investment costs $AC^{INV,DER}$ and operational costs $C_s^{OM,DER}$ associated with DERs, electricity purchasing costs $C_s^{grid}$, and income from generating and/or selling renewable energy $I_s$:

$$\min \quad TAC = \sum_{DER} AC^{INV,DER} + \sum_{s \in S} \left( C_s^{grid} + \left( \sum_{DER} C_s^{OM,DER} \right) + I_s \right) \quad (5)$$

This study considers PVs, batteries, ASHPs, and gas boilers as DERs. At each dwelling $i \in I$ and time point $t \in T$, the base electricity demand $E_{i,t}^{load}$ and the consumption by ASHPs $E_{i,t,p}^{hp}$ are met by electricity purchased from the grid $E_{i,t}^{grid}$, PV generation $E_{i,t}^{PV}$, and/or battery discharge $E_{i,t,c}^{discharge}$:

$$\hat{E}_{i,t}^{load} + \sum_p E_{i,t,p}^{hp} = E_{i,t}^{grid} + E_{i,t}^{PV} + \sum_c E_{i,t,c}^{discharge} \quad (6)$$

Note that each ASHP is denoted by subscript $p \in P$, while batteries are denoted by subscript $c \in C$. All parameters are non-italicised and emphasised using the circumflex symbol.

A logical formulation is enforced in De Mel et al. [13] to prevent the simultaneous buying and selling of electricity. This is formulated using binary variable $X_{i,t}$ in the MILP formulation. In the presence of ASHPs, the following modified big-M constraints are utilised to avoid nonlinearity.

$$E_{i,t}^{grid} \leq \left( \widehat{E}_{i,t}^{load} + \sum_p E_{i,t,p}^{hp} \right) \tag{7}$$

$$E_{i,t}^{grid} \leq M \cdot (1 - X_{i,t}) \tag{8}$$

$$E_{i,t}^{sold} \leq M \cdot X_{i,t} \tag{9}$$

Thermal storage, usually in the form of hot water cylinders or tanks, is usually required where ASHPs meet both space heating and hot water demands [47]. The heat demand $\widehat{H}_{i,t}^{load}$ can be met via thermal discharge from hot water tanks $H_{i,t,k}^{disch}$ and heat generation from boilers $H_{i,t,b}$:

$$\widehat{H}_{i,t}^{load} = \sum_k H_{i,t,k}^{disch} + \sum_b H_{i,t,b} \tag{10}$$

Note that subscripts $k \in K$ and $b \in B$ denote the discrete technology selections of hot water tanks and gas boilers. All heat generated by the ASHPs is used to charge the respective hot water tank $H_{i,t,k}^{ch}$, as shown below.

$$H_{i,t,k}^{ch} = \sum_p H_{i,t,p,k}^{hp} \tag{11}$$

The ASHP-tank combination is described by a single binary variable $J_{i,p,k}$, where one ASHP-tank combination is permitted for installation at each dwelling:

$$\sum_{k,p} J_{i,p,k} \leq 1 \tag{12}$$

The use of one binary variable ensures that tanks that are incompatible with the rated heat pump output are not chosen, as specific requirements such as flowrates cannot be met. This is further enabled using a parametric indicator $I_{p,k} \in \{0,1\}$, which indicates whether the heat pump is feasible for operation with the tank as indicated in manufacturer datasheets:

$$J_{i,p,k} \leq I_{p,k} \tag{13}$$

The electricity consumption of ASHPs $E_{i,t,p}^{hp}$ with respect to heat generation can be described using the Coefficient of Performance $\widehat{COP}_{t,p}$:

$$E_{i,t,p}^{hp} = \frac{\sum_k H_{i,t,p,k}^{hp}}{\widehat{COP}_{t,p}} \tag{14}$$

The heat generation capability of an ASHP is dependent on ambient temperature. As a result, the maximum heat generation of an ASHP $H_{t,p}^{hp,max}$ is restricted by ambient temperature:

$$\sum_k H_{i,t,p,k}^{hp} \leq \widehat{H}_{t,p}^{hp,max} \cdot \sum_k J_{i,k,p} \tag{15}$$

Similarly, a minimum operating temperature $T_p^{min}$ is often specified for each ASHP, below which operation ceases. This is described by the constraint below:

$$If \ \widehat{T}_t^{air} \leq \widehat{T}_p^{min} \quad \sum_k H_{i,t,p,k}^{hp} = 0 \tag{16}$$

The determination of $\widehat{COP}_{t,p}$ and maximum capacity with respect to varying ambient temperatures is further detailed in Appendix 6A.

The annualised investment cost associated with ASHPs $C^{inv,hp}$ is calculated as shown below, which includes technology unit cost $\widehat{CU}_p$ and installation cost $\widehat{CI}_p$.

$$C^{inv,hp} = CRF \cdot \sum_{i,p} \left(\widehat{CU}_p + \widehat{CI}_p\right) \cdot \sum_k J_{i,p,k} \tag{17}$$

Note that CRF denotes the capital recovery factor, which has been used to annualise the costs.

Focusing on thermal storage, the heat balance of a non-stratified tank can be described as below, based on the formulation by Sarbu and Sebarchievici [48]:

$$If\ t > 1 \quad \frac{\hat{c}_p \widehat{V}_k \hat{\rho}\left(T_{i,t,k}^{tank} - T_{i,t-1,k}^{tank}\right)}{\Delta t * 3600} = H_{i,t,k}^{ch} \cdot \hat{\eta}_k^{ch} - \frac{H_{i,t,k}^{disch}}{\hat{\eta}_k^{disch}} - H_{i,t,k}^{loss}$$

(18)

$$else \quad \frac{\hat{c}_p \widehat{V}_k \hat{\rho}\left(T_{i,t,k}^{tank} - \widehat{T}_t^{air}\right)}{\Delta t * 3600} = H_{i,t,k}^{ch} \cdot \hat{\eta}_k^{ch} - \frac{H_{i,t,k}^{disch}}{\hat{\eta}_k^{disch}} - H_{i,t,k}^{loss}$$

where $\hat{c}_p$ is the specific heat capacity, $T_{i,t,k}^{tank}$ is the tank temperature, $\widehat{V}_k$ is the volume of the tank, $\hat{\rho}$ is the density of water, and $\hat{\eta}_k^{ch}$ and $\hat{\eta}_k^{ch}$ are charging and discharging efficiencies. The variable $H_{i,t,k}^{loss}$ represents the heat loss from the tank, which is described as a function of the average heat loss $\widehat{H}_k^{loss}$ recorded for each tank in the manufacturer's datasheets.

$$H_{i,t,k}^{loss} = \widehat{H}_k^{loss} \cdot \sum_p J_{i,p,k} \quad (19)$$

The temperature in the tank is specified an upper and lower bound, where the upper bound is the pre-specified heat pump supply temperature $\widehat{T}^{ws}$, and the lower bound is $\widehat{T}_k^{min}$:

$$T_{i,t,k}^{tank} \leq \widehat{T}^{ws} \cdot \sum_p J_{i,p,k} \quad (20)$$

$$T_{i,t,k}^{tank} \leq \widehat{T}_k^{min} \cdot \sum_p J_{i,p,k} \quad (21)$$

The temperature of the tank at the start and end of the overall time period considered for each representative day must be the same, to ensure consistency across all the days of a season:

$$T_{i,t=start,k}^{tank} = T_{i,t=end,k}^{tank} \quad (22)$$

Further logical big-M constraints are imposed to ensure that tanks that are not selected do not charge nor discharge heat:

$$H_{i,t,k}^{ch} \leq \mathrm{M} \cdot \sum_p J_{i,p,k} \tag{23}$$

$$H_{i,t,k}^{disch} \leq \mathrm{M} \cdot \sum_p J_{i,p,k} \tag{24}$$

The annualised investment cost of the tanks can be described by the equation below, which accounts for the unit cost $\widehat{\mathrm{CU}}_{p,k}$ of the tank:

$$C^{inv,tank} = \mathrm{CRF} \cdot \sum_{i,k} \sum_p J_{i,p,k} \cdot \widehat{\mathrm{CU}}_{p,k} \tag{25}$$

Note that, in this case, the unit cost of the tank varies with the ASHP chosen, potentially due to variations in fittings and pipework. This is further described in Appendix A. It is assumed that the installation cost of the tank is accounted for in the installation cost of the ASHP.

With respect to the battery formulation suitable for discrete capacities, the battery unit (with a respective capacity) chosen is restricted to one at each dwelling using a binary variable $W_{i,c}$:

$$\sum_c W_{i,c} \leq 1 \tag{26}$$

Battery charging and discharging is restricted to the maximum charging and discharging power $\mathrm{E}_c^{max}$, as provided in manufacturer datasheets:

$$E_{i,t,c}^{ch} \leq \mathrm{E}_c^{max} * W_{i,c} \tag{27}$$

$$E_{i,t,c}^{disch} \leq \mathrm{E}_c^{max} \cdot W_{i,c} \tag{28}$$

Other constraints include ensuring that the energy stored in the battery do not exceed state of charge or depth of discharge. The annualised investment costs are described by unit cost $\widehat{\mathrm{CU}}_c$ and installation cost $\widehat{\mathrm{CI}}_c$:

$$C^{INV,batt} = \mathrm{CRF} \cdot \sum_{i,c} W_{i,c} \left( \widehat{\mathrm{CU}}_c + \widehat{\mathrm{CI}}_c \right) \tag{29}$$

The operational cost (per season) considers a fixed operational cost per year $CO_c$ into account, where $N^{days}$ represents the number of days in a season:

$$C^{OM,batt} = \sum_{i,c} W_{i,c} \cdot \widehat{CO}_c \cdot \frac{1}{365} \cdot N^{days} \tag{30}$$

The discrete boiler formulation includes a binary variable $U_{i,b}$, which indicates the unit (and respective capacity) chosen. Boiler selection is also restricted to a maximum of 1 per household:

$$\sum_b U_{i,b} \leq 1 \tag{31}$$

Heat generated by the boiler $H_{i,t,b}$ is restricted to its upper bound, a specified maximum capacity $\widehat{H}_b^{max}$:

$$H_{i,t,b} \leq \widehat{H}_b^{max} * U_{i,b} \tag{32}$$

The annualised investment cost of boilers is calculated below, where $CU_b$ is the unit cost and $CI_b$ is the installation cost:

$$C^{INV,B} = \text{CRF} \cdot \sum_{i,b} U_{i,b} \left(\widehat{CU}_b + \widehat{CI}_b\right) \tag{33}$$

The operational cost of the boiler accounts for gas usage through the gas price $\widehat{C}^{gas}$ (£/kWh) and boiler efficiency $\hat{\eta}_b$. Note that $\Delta t$ represents the discretised time interval.

$$C^{OM,B} = \sum_{i,t,b} H_{i,t,b} * \Delta t * \frac{\widehat{C}^{gas}}{\hat{\eta}_b} \cdot N^{days} \tag{34}$$

The reader is directed to previous work [13] for other general constraints, continuous PV formulation, complementarity constraints, nonconvex MOPF formulation, and linking constraints.

## 3. Computational Results and Discussion

The overall algorithm for obtaining discrete DES designs with respect to MOPF is tested with complementarity reformulations algorithms $CR$ and $CR$-$H$. The performance of these algorithms is also compared to the existing commercial MINLP Branch-and-Bound solver SBB [49]. As SBB is a generic solver and has not been customised for DES-MOPF applications (unlike the above algorithms), it is provided a feasible mixed-integer initial solution by solving the MILP first, followed by the NLP with fixed binary variables and MOPF constraints, to assist with solution of the overall optimisation problem. The following notation is used to describe the solution methods to aid comparison and discussion:

**PA:** The proposed overall algorithm with the original complementarity reformulations algorithm $CR$. This solves the combined DES-MOPF design formulation.

**PA-H:** The overall algorithm with the complementarity reformulations algorithm $CR$-$H$, which includes a heuristic modification that terminates the complementarity algorithm when poor relaxed objective functions are found.

**SBB**: Commercial MINLP solver SBB [49] with a feasible mixed-integer initialisation.

Additionally, **MILP** refers to the DES design model which is used for initialising the algorithms and determining the binary decision variables. This includes linear DC approximations in place of MOPF constraints. The maximum time limit for both PA and PA-H is set to 10 hours, excluding initialisations. These algorithms use commercial solvers CPLEX [50] and CONOPT [51] for the decomposed MILP and NLP problems, respectively. The SBB solver time limit is also set to 10 hours, excluding the time taken for initialisations, while leaving other solver options as default. All computational experiments were carried out on an Intel® Core™ i7-10510U CPU at 1.80GHz - 2.30 GHz. All code and input files are provided in an external repository on Github: https://github.com/Ishanki/Discrete-DES-MOPF

### 3.1 Case Studies

The solution methods are tested using DES that are located within low-voltage (LV) and unbalanced radial distribution networks, adapted from the 906-node IEEE European LV test feeder [52] with a -30° phase shifting Delta-Wye ($\Delta - Y$) transformer. Network 1, which includes 96 buses in total out of which 12 are consumer or load nodes, is presented in Figure

4. Network 2, which is an extension of Network 1 with 256 buses out of which 22 are consumer or load nodes, is presented in Figure 5. Note that each of the consumers are single-phase loads, connected to the three-phase LV network using single-phase feeders, which is typical of UK LV networks. Peak electricity and heat demands at each of these consumer nodes, as found in the winter season, are provided in Table 1. Averaged 24-hour seasonal demand profiles for electricity and heating have been generated for each of the of the consumers. Solar irradiance and ambient temperature profiles have been obtained for Surrey, UK [53], which are presented in Figure 6. All demand and weather profiles are discretised into hourly timestamps. PVs are retained as continuous decisions, as these are dependent on the area. The maximum total PV area that can be installed by each consumer is set to 35 m$^2$. Commercially available types and capacities for lithium-ion batteries, boilers, and ASHPs are included as discrete technologies. Parameters for the discrete technologies used are provided in Appendix A, while scalar parameters for pricing and PVs are provided in Appendix B. Note that both a generation tariff and export tariff are considered in this case study, to encourage the installation of solar panels and enable the testing of the algorithms.

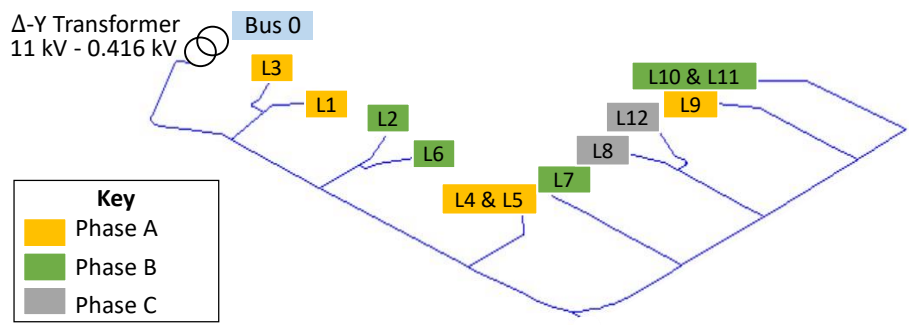

*Figure 4. Network 1, which includes 12 load/consumer nodes.*

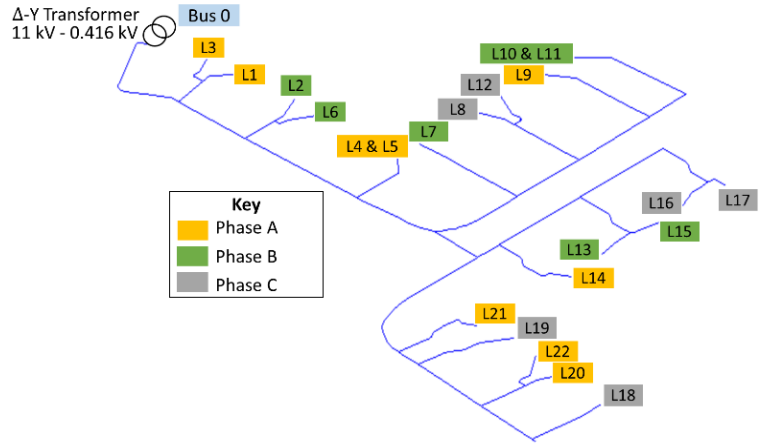

*Figure 5. Network 2, with 22 load/consumer nodes*

*Table 1. Peak electricity and heat demands, as found in the winter season, at consumer nodes.*

| Load | Peak Electricity Demand (kW) | Peak Heat Demand (kW) |
|---|---|---|
| L1 | 1.39 | 12.2 |
| L2 | 2.7 | 23.7 |
| L3 | 1.2 | 10.5 |
| L4 | 1.58 | 14 |
| L5 | 1.08 | 9.5 |
| L6 | 0.89 | 7.8 |
| L7 | 1.88 | 16.6 |
| L8 | 1.83 | 16.1 |
| L9 | 2.86 | 25.2 |
| L10 | 2.44 | 21.5 |
| L11 | 1.08 | 9.5 |
| L12 | 0.51 | 4.5 |
| L13 | 2.17 | 19.1 |
| L14 | 0.63 | 5.6 |
| L15 | 1.01 | 8.9 |
| L16 | 1.28 | 11.3 |
| L17 | 0.72 | 6.4 |
| L18 | 1.85 | 16.3 |
| L19 | 1.52 | 13.4 |
| L20 | 1.97 | 17.3 |
| L21 | 1.23 | 10.8 |
| L22 | 0.6 | 5.2 |

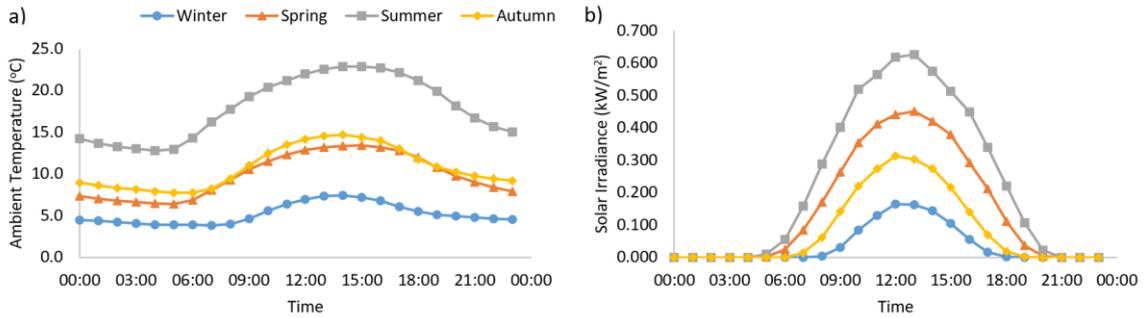

*Figure 6. Weather profiles used in the case studies: a) ambient temperature and b) solar irradiance.*

The following case studies are used to test the proposed solution methods, for which results are presented in the subsequent subsections:

**Case 1a:** Network 1 with the choice of PVs, batteries, ASHPs, and hot water (HW) tanks.

**Case 1b:** Network 1 with the choice of PVs, batteries, and gas boilers.

**Case 1c:** Network 1 with the choice of PVs, batteries, boilers, ASHPs, and HW tanks. An exceptionally high renewable electricity generation tariff is used to further incentivise local PV generation.

Note that ASHPs and boilers are provided as separate choices in Case 1a and 1b because, under a cost objective function, higher investment cost of ASHPs and tanks result in boilers always being selected, as made evident in Case 1c. The MILP solved within the algorithms, which contains DES design constraints for all DERs and the linear DC approximations for power flow, includes 85,252 continuous variables, 6,720 binary variables, and 123,945 inequality and equality constraints. The NLP with all DERs, which includes MOPF constraints and reformulated complementarity constraints, includes 204,869 continuous variables and 246,381 inequality and equality constraints.

Cases 2a – c, which involve the larger network, Network 2, are used to further test the efficacy of PA-H.

**Case 2a:** Network 2 with the choice of PVs, batteries, ASHPs, and HW tanks.

**Case 2b:** Network 2 with the choice of PVs, batteries, and gas boilers.

**Case 2c:** Network 2 with ASHPs and HW tanks only (no local electricity generation).

The MILP with all DERs includes 171,693 continuous variables, 12,320 binary variables, and 258,145 inequality and equality constraints. The NLP with all DERs includes 490,989 continuous variables and 590,291 inequality and equality constraints.

## 3.2 Performance of solution methods

The results for the case studies using Network 1 (Cases 1a – c) using each solution method are presented in Table 2. The initial MILP objective value (without MOPF constraints), final DES-MOPF objective value from each solution method, percentage difference of the final objective values with respect to the initial MILP value, total CPU time taken (including initialisations) and the final solve status are reported here.

For all three cases, where locally optimal integer solutions have been found, the percentage difference in the DES-MOPF objective value when compared to that of the MILP (which includes DES design with linear DC approximations) is negligible. This may suggest that MOPF constraints are unimportant in DES design by observing the objective values alone, however, the impacts on the designs are further investigated in Section 3.3 below.

The proposed algorithms outperform the commercial MINLP solver SBB. In Case 1a and 1c, locally optimal and mixed-integer solutions are found by proposed algorithms PA and PA-H where SBB fails to find any. Note that SBB terminates in Case 1a and Case 1c due to failing to find locally optimal solutions at the root node, where a relaxed MINLP is solved. For Case 1b, a slightly better objective value is obtained by PA and PA-H within ~7-12 times less computational time.

*Table 2. Performance of the solution methods tested under different test cases. Note that the percentage difference of objective values is calculated with respect to the initial MILP objective value.*

| Case | Algorithm | Initial MILP Objective (£) | Final Objective (£) | % Difference in obj. values | Total CPU Time taken (s) | Final Algorithm Status |
|---|---|---|---|---|---|---|
| Case 1a | PA | 46,431 | 46,488 | 0.12 | 9,112 | Converged |
|  | PA-H | 46,431 | 46,488 | 0.12 | 2,727 | Converged |
|  | SBB | 46,431 | - | - | 1,990 | No solution |
| Case 1b | PA | 13,492 | 13,572 | 0.60 | 2,908 | Converged |
|  | PA-H | 13,492 | 13,572 | 0.60 | 1,866 | Converged |
|  | SBB | 13,492 | 13,588 | 0.71 | 23,137 | Not converged – default node limit exceeded |
| Case 1c | PA | -145,571 | -144,279 | 0.89 | 39,106 | Not converged - time limit exceeded |
|  | PA-H | -145,571 | -144,542 | 0.71 | 37,852 | Not converged - time limit exceeded |
|  | SBB | -145,571 | - | - | 30,781 | No solution |

In Cases 1a and 1b, the proposed algorithm PA and its heuristic modification PA-H both converged within a few iterations (i.e., the lower bound exceeded the lowest upper bound *LUB*) and produced locally optimal solutions. The progression of the upper and lower bounds with respect to each proposed algorithm and iteration is portrayed in Figure 7. In both these cases, the designs proposed in the first iteration produced the lowest upper bounds *LUB* with respect to the NLP.

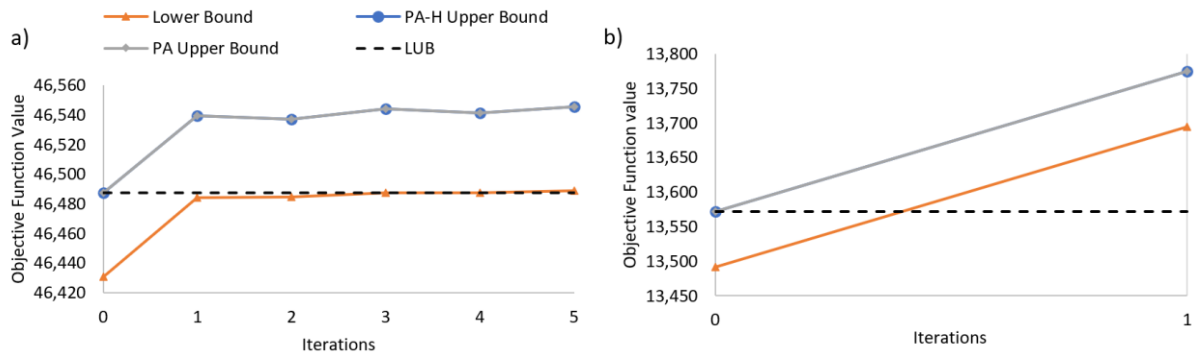

*Figure 7. Lower and upper bounds for PA and PA-H algorithms for a) Case 1a and b) Case 1b.*

The results of Case 1c are different to the above two cases as neither algorithm converges, albeit finding feasible and locally optimal solutions. The progression of the upper and lower bounds for Case 1c are portrayed in Figure 8, and the *LUB* for each algorithm is indicated as black and red points. Note that *LUB*s for the two proposed algorithms are different, unlike in the previous cases. This is because the lower computational expense of PA-H led to a higher total number of iterations than PA before the maximum time limit was reached. The lower bound appears to increase relatively slowly, and in some instances, there is very little increase in objective value (<0.0001%) despite the integer cuts ensuring that the designs change at each iteration. This indicates the presence of symmetrical solutions with respect to the MILP. The instability of the upper bound can be clearly observed in Figure 8, where sudden decreases in the upper bound do not necessarily correspond to step changes in the lower bound. Note that these fluctuations in the upper bound do not necessarily reflect the presence of local optima as changes in the binary topology are also reflected by the upper bound. The design changes that cause significant improvements to the upper bound are further investigated in Section 3.3 below.

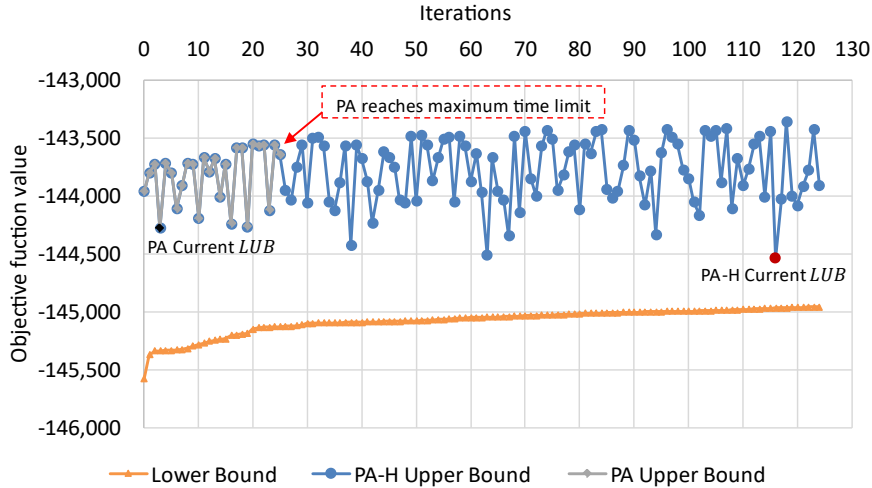

*Figure 8. The behaviour of the lower bound LB and upper bound UB for Case 1c. The black dot represents the current LUB for algorithm PA, and the red dot is the same for algorithm PA-H.*

The efficacy of the heuristic modification in PA-H is assessed by evaluating the upper bounds produced at each iteration in comparison with PA. The solutions are similar at each iteration between the two algorithms, and this is made evident in Figure 7 and Figure 8. In Case 1c (Figure 8), the lower computational expense of PA-H led it to complete many more iterations compared to PA prior to reaching the maximum time limit and terminating. This suggests that PA-H is a good alternative to PA when a balance between accuracy and computational expense is required.

The efficacy of PA-H on the larger Network 2 is further tested using Cases 2a – c, for which results are presented in Table 3. Unlike in Cases 1a and 1b, Cases 2a and 2b do not converge within the time limit. This is primarily because these cases have a greater number of variables and constraints compared to Cases 1a and 1b, which causes each NLP solve to take approximately 5 times longer on average. A larger gap between the objectives of the initial MILP prediction and the DES-MOPF is also observed. However, feasible upper bounds are available for Cases 2a and 2b. Case 2c, which only includes ASHPs and HW tanks, converges and returns the same objective value as the initial MILP.

*Table 3. Performance of PA-H under Cases 2a, b, and c.*

| Case | MILP Objective (£) | Final objective value (£) | % Difference in obj. values | Total CPU Time taken (s) | Status |
|---|---|---|---|---|---|
| Case 2a | 79,246 | 80,281 | 1.31 | 44,188 | Not converged - time limit exceeded |
| Case 2b | 22,813 | 23,688 | 3.83 | 44,340 | Not converged - time limit exceeded |
| Case 2c | 88,850 | 88,850 | 0.00 | 7,453 | Converged |

Figure 9 portrays the progression of the lower and upper bounds for Cases 2a and 2b. Note that Case 2c converges prior to requiring any integer cuts and is not presented graphically, as the ASHP-tank selections from the MILP are feasible with respect to the network constraints as well. The existing *LUB* for both Cases 2a and 2b have been obtained at Iteration 0, and no improved upper bounds have been obtained in subsequent iterations. The upper bound in Case 2a (Figure 9a) is observed to be generally very stable, with step changes in the lower bound (due to changes in the binary topology) reflected clearly on the upper bound as well. Notice the increase in the lower bound is very slow, once again indicating the presence of symmetrical solutions. On the other hand, the upper bound for Case 2b (Figure 9b) is less stable, suggesting that changes in the binary topology have more impact on the MOPF constraints.

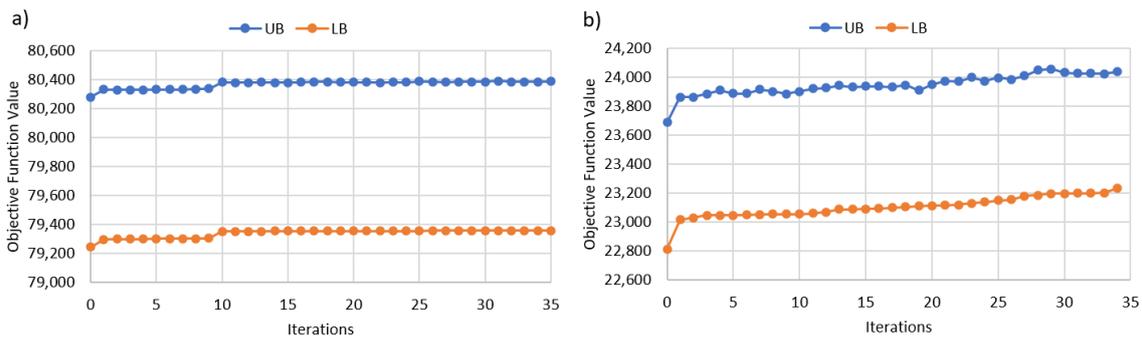

*Figure 9. Lower and upper bounds for PA-H obtained for a) Case 2a and b) Case 2b.*

## 3.3 Impacts on DES design

Table 4 presents the annualised costs for Cases 1a and 1b, which includes all investment costs, operational costs, and income generated by the different solution methods tested in the study. These costs reflect the design decisions made. As previously discussed, the solutions for the proposed algorithms PA and PA-H are identical for these two cases.

In both cases, no batteries have been selected due to high investment and operational costs. A very minor decrease in PV investment is observed in the results of PA and PA-H, when compared to the initial MILP prediction where the maximum possible PV area is selected (which does not include MOPF constraints). This is reflective of the reduced total PV area in the results proposed by PA and PA-H. Generation and export income are also slightly reduced. These reductions are observed because the MOPF constraints are reducing the excess PV power exports to the network to eliminate violations to the voltage upper bound imposed by the network. This is further confirmed by a post-optimisation power flow calculation (with all key DES design and operational variables fixed), which returns an infeasible solution for the combined design and operational schedules recommended by the initial MILP.

*Table 4. Annualised costs for Case 1a and 1b.*

| Annualised Costs (£) | Case 1a | | | | Case 1b | | | |
|---|---|---|---|---|---|---|---|---|
| | MILP | PA | PA-H | SBB | MILP | PA | PA-H | SBB |
| Objective value | 46,431 | 46,488 | 46,488 | - | 13,492 | 13,572 | 13,572 | 13,588 |
| Electricity purchase | 23,943 | 23,953 | 23,953 | - | 3,029 | 3,047 | 3,047 | 3,043 |
| PV investment | 10,594 | 10,534 | 10,534 | - | 10,594 | 10,164 | 10,164 | 10,378 |
| PV operation | 750 | 746 | 746 | - | 750 | 720 | 720 | 735 |
| Boiler investment | 0 | 0 | 0 | - | 2,827 | 2,827 | 2,827 | 2,827 |
| Boiler operation | 0 | 0 | 0 | - | 8,097 | 8,097 | 8,097 | 8,097 |
| Battery investment | 0 | 0 | 0 | - | 0 | 0 | 0 | 0 |
| Battery operation | 0 | 0 | 0 | - | 0 | 0 | 0 | 0 |
| ASHP investment | 14,591 | 14,591 | 14,591 | - | 0 | 0 | 0 | 0 |
| HW Tank investment | 7,438 | 7,438 | 7,438 | - | 0 | 0 | 0 | 0 |
| Export income | -2,513 | -2,477 | -2,477 | - | -3,433 | -3,261 | -3,261 | -3,331 |
| Generation income | -8,372 | -8,296 | -8,296 | - | -8,372 | -8,020 | -8,020 | -8,160 |

Interestingly, Case 1a with ASHPs instead of gas boilers supports nearly 4% higher PV integration compared to Case 1b (comparing the results of PA across both case studies), due to the fact that more energy generated by the PVs is consumed locally by the ASHPs at each consumer node. This is also reflected in the differences in export income, where Case 1a exports nearly 28% less electricity compared to Case 1b with respect to the results from the proposed algorithms. This suggests that PV-ASHP installation is more advantageous to the underlying unbalanced LV network compared to PV-Boiler installations, as the higher proportion of local generation and consumption minimises the power export and consequent network violations. Unfortunately, despite the potential benefits to the network, this results in much higher investment and operational costs to the consumers in Case 1a.

The annualised costs of Case 1c are shared in Table 5. Note that the results for the current lowest upper bounds (*LUB*s) are reported for PA and PA-H, as the algorithms did not converge. As mentioned previously, an extremely high renewable energy generation tariff was used in this test case. This has led to the highest possible PV area being recommended by the MILP, PA, and PA-H, resulting in the same PV investment costs as well. However, unlike in the previous cases, both PA and PA-H have chosen to install batteries in this case to minimise power exports and subsequent voltage violations. The integer cuts present in the algorithms have varied battery capacities and their locations at each iteration, leading to the fluctuations observed in the upper bound in Figure 8 of Section 3.2 above. Installing batteries has contributed to finding better upper bounds, as batteries can reduce the power exported to the network by storing excess energy at peak PV production times and discharging for local consumption later in the day. However, the fluctuations show that not all battery installations produce the same effect on the upper bound. This is made evident in Figure 10, where battery locations and sizes at loads L1-L12 are shown along with the PA upper bound at each iteration. Looking at this figure, it is not possible to predetermine the optimal locations and battery sizes for the DES to maximise PV generation while eliminating voltage violations in the LV network. This highlights the importance of using an optimisation framework or algorithm, such as PA or PA-H, capable of integrating nonconvex MOPF constraints to ensure that feasible DES designs are obtained.

*Table 5. Annualised costs for Case 1c.*

| Annualised Costs (£) | Case 1c | | | |
|---|---|---|---|---|
| | MILP | PA | PA-H | SBB |
| Objective value | -145,571 | -144,279 | -144,542 | - |
| Electricity purchase | 3,029 | 2,866 | 2,656 | - |
| PV investment | 10,594 | 10,594 | 10,594 | - |
| PV operation | 750 | 750 | 750 | - |
| Boiler investment | 2,827 | 2,827 | 2,827 | - |
| Boiler operation | 8,097 | 8,097 | 8,097 | - |
| Battery investment | 0 | 248 | 609 | - |
| Battery operation | 0 | 110 | 270 | - |
| ASHP investment | 0 | 0 | 0 | - |
| HW Tank investment | 0 | 0 | 0 | - |
| Export income | -3,433 | -3,361 | -3,314 | - |
| Generation income | -167,435 | -166,409 | -167,030 | - |

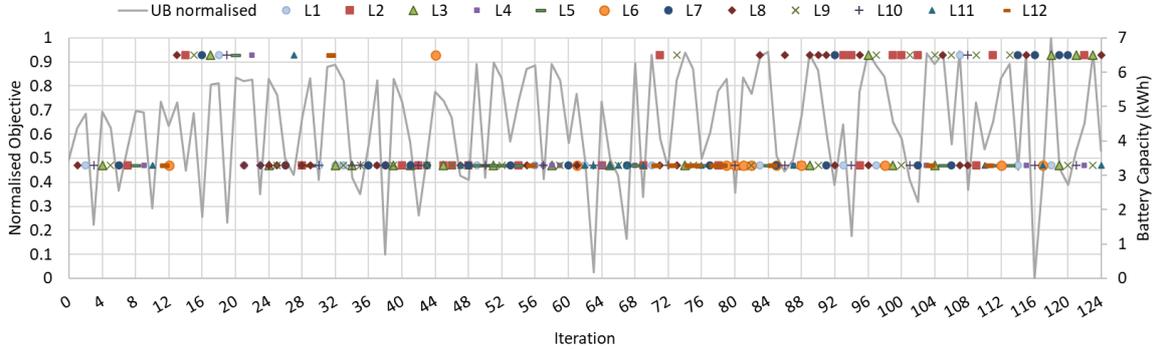

*Figure 10. Locations of batteries and their respective sizes, alongside the normalised upper bound obtained from the solution method PA.*

Table 6 presents the annualised costs of Cases 2a – c for Network 2, which have been obtained via the initial MILP and PA-H algorithm. Note that for Cases 2a and 2b, the current *LUB*s have been reported as the algorithm did not converge for either of these cases. The trends between Cases 2a and 2b are similar to those observed for Cases 1a and 1b, where generation and export income has been reduced compared to the initial MILP prediction to eliminate voltage violations. The PA-H solution for Case 2a, which includes ASHPs and HW tanks, supports 16% higher PV capacity compared to that of Case 2b with boilers, as more locally generated electricity is consumed by the ASHPs. However, Case 2a is over three times more expensive than Case 2b, primarily due to the high investment and operational costs of ASHPs considered in this case study.

*Table 6. Annualised costs for Cases 2a, 2b, and 2c, as obtained using algorithm PA-H.*

| Annualised Costs (£) | Case 2a | | Case 2b | | Case 2c | |
| --- | --- | --- | --- | --- | --- | --- |
| | MILP | PA-H | MILP | PA-H | MILP | PA-H |
| Objective value | 79,246 | 80,281 | 22,813 | 23,688 | 88,850 | 88,850 |
| Electricity purchase | 39,799 | 40,215 | 5,089 | 5,324 | 50,084 | 50,084 |
| PV investment | 19,422 | 17,165 | 19,422 | 14,734 | 0 | 0 |
| PV operation | 1,375 | 1,215 | 1,375 | 1,043 | 0 | 0 |
| Boiler investment | 0 | 0 | 5,152 | 5,152 | 0 | 0 |
| Boiler operation | 0 | 0 | 13,526 | 13,526 | 0 | 0 |
| Battery investment | 0 | 0 | 0 | 0 | 0 | 0 |
| Battery operation | 0 | 0 | 0 | 0 | 0 | 0 |
| ASHP investment | 25,129 | 25,129 | 0 | 0 | 25,129 | 25,129 |
| HW Tank investment | 13,637 | 13,637 | 0 | 0 | 13,637 | 13,637 |
| Export income | 4,768 | 3,831 | 6,403 | 4,554 | 0 | 0 |
| Generation income | 15,348 | 13,249 | 15,348 | 11,538 | 0 | 0 |

For Case 2c, which relies solely on purchasing electricity from the grid to power ASHPs, the initial MILP design is also feasible with respect to the DES-MOPF framework. Case 2c produces the highest objective value out of the three cases. This is because of the lack of generation and export income from PV production, in contrast to Case 2a. Note that the generous tariffs included in these case studies for local electricity export and generation help offset PV investment and operational costs for the consumers in Case 2a, while reducing electricity purchasing costs.

In general, these results imply that greater local renewable energy consumption is preferred over export to the grid to help minimise the impacts to unbalanced distribution networks. This suggests that incentivising the installation and/or operation of electrified heating systems (such as ASHPs) instead of electricity export for consumers with DERs is more beneficial to network operators as well.

Note that this study has focused on testing and demonstrating the framework and algorithms, as opposed to providing specific recommendations for policy-making. Doing the latter will require assessing the uncertainty of the parameters used, such as the highly volatile gas prices observed in 2022, which can alter design decisions and cost predictions. This is especially relevant for technology comparisons, as scenarios with ASHPs and HW tanks are over three times more expensive than those with boilers, based on the parameters used here. These results do, however, highlight that significant policy considerations and incentives will be required to not only help minimise the impacts on the unbalanced distribution networks, but also facilitate the decarbonisation of fossil fuel-based residential heating systems.

### 3.4 Limitations of the proposed methods

The proposed solution methods enable the inclusion of discrete technologies such as ASHPs and batteries, to obtain locally optimal and feasible solutions with respect to nonconvex MOPF constraints representing the underlying LV unbalanced networks. Results demonstrate the importance of doing so, particularly in Case 1c where discrete batteries have a significant impact on the upper bound produced by the NLP. However, these methods can be further improved to ensure that they can be tested against a wide range of networks and case studies.

Although the proposed methods perform significantly better than the existing MINLP commercial solver by finding feasible solutions, larger case studies require more computational time and expense before achieving convergence. PA-H is a much more efficient alternative to PA, as it produces the same results while taking up to 70% less computational time for Case 1 studies. However, it too fails to converge for Case 1c before all feasible binary topologies are explored, suggesting that there may be better upper bounds in the unexplored binary topologies. The lack of convergence is also observed for larger Cases 2a and 2b which consider Network 2 with over twice more variables and constraints. This can hinder the scalability of the proposed methods, as case studies with larger number of buses or discrete options may take much longer to converge. To improve the scalability and computational efficiency, other methods for obtaining tighter lower bounds need to be explored, or additional decompositions and further heuristic modifications will be required. Such improvements can also enable uncertainty analyses to assess which parameters have the most influence of DES designs, while considering MOPF constraints.

The DES-MOPF formulations in this study use nonconvex MOPF constraints. This is primarily because using nonconvex power flow constraints can guarantee accuracy with respect to network constraints, as opposed to using linear or convex approximations. Furthermore, as previously mentioned, existing multiphase linearisations require prior knowledge of the design for feasible power flow initialisations, and/or feasible solutions from external power flow simulation tools at each iteration. The proposed methods eliminate any potentially expensive computational calls to external tools or programmes, as feasible initialisations are embedded in the algorithms. However, further research is required to assess if such iterative linearisations are more suitable alternatives to using integer cuts and complementarity reformulations, in terms of solution accuracy and computational speed.

## 4. Conclusions

Existing optimisation frameworks for the design of DES have often overlooked the role of the unbalanced distribution networks to which they are connected. Thus far, no optimisation framework has been proposed to find discrete designs for DES, which are also feasible with multiphase network constraints. This study proposes new optimisation methods to obtain discrete designs for DES, while including MOPF constraints representing the underlying unbalanced distribution network. The proposed solution methods focus on obtaining feasible and locally optimal solutions by decomposing the overall MINLP formulation to MILP and Nonlinear Programming (NLP) subproblems, where integer cuts, complementarity reformulations, and existing MILP and NLP commercial solvers are utilised. A heuristic modification to the first algorithm is proposed to minimise computational time. Conventional distributed technologies such as PVs, batteries, and gas boilers are included in the formulations. Improved formulations for electricity-consuming ASHPs that operate in conjunction with hot water storage are also included, as an alternative to traditional gas boilers. The algorithms are tested using two low-voltage unbalanced networks of varying size under six different cases (3 for each network).

Results for three cases using the smaller network are compared with a commercial MINLP solver, SBB, and the proposed algorithms outperform SBB in every case. Locally optimal solutions are found for all cases, while SBB fails to find any solutions for two out of the three cases. Where convergence criteria are met, the algorithm with the metaheuristic modification improves solution time by up to 70% in comparison with the algorithm without this modification. In one case where an exceptionally high renewable generation tariff is set, both proposed algorithms fail to converge before exceeding the maximum time limit, but find feasible solutions. This is primarily due to the influence of discrete battery decisions, where results show that the location and sizes of the battery can play a huge role in whether the overall design is compatible with network constraints or not. The proposed algorithm with the heuristic modification is further tested using three cases for a larger network. While feasible upper bounds are obtained, two out of the three cases do not achieve convergence due to larger computational expense and a larger gap observed between the upper and lower bounds. For the parameters used in this study, results for both networks also suggest that the combination of ASHPs with PVs is more conducive to increasing local renewable energy generation while respecting network constraints, than the less-expensive combination of gas boilers and PVs. These results also further demonstrate the necessity of considering MOPF constraints when designing DES, as not doing so could lead to greater costs for the consumer and reduced renewable generation. The scalability and convergence of proposed methods can be improved in future work by incorporating additional decompositions and/or heuristics. The formulations can also be extended to include additional features related to future smart grids, such as capacitors and network switches, or nonlinear phenomena related to DES technologies, such as battery degradation.

## Appendix A – Parameters of discrete technologies

All ASHPs presented in Table A.7 are Mitsubishi Ecodan heat pumps using refrigerant R32 [54]. Note that charging and discharging efficiencies ($\hat{\eta}_k^{ch}$ and $\hat{\eta}_k^{disch}$) have been assumed, and $\hat{T}_k^{min}$ is set as 10°C lower than the preset ASHP water supply temperature, $\hat{T}^{ws} = 55°C$.

*Table A.7. Parameters for the discrete selection of ASHPs.*

| Label | Model | $\hat{T}_p^{min}$ (°C) | $\widehat{CU}_p$ (£) | $\widehat{CI}_p$ (£)* | References |
|---|---|---|---|---|---|
| S1 | PUZ-WM50VHA(-BS) | -10 | 2,575 | 1,931 | [54], [55] |
| M1 | PUZ-WM60VAA(-BS) | -10 | 3,053 | 2,137 | [54], [55] |
| M2 | PUZ-HWM140VHA(-BS) | -15 | 4,861 | 3,403 | [54], [55] |
| L1 | 2 x PUZ-HWM140VHA(-BS) | -15 | 9,722 | 6,805 | [54], [55] |

*Installation costs assumed to be 70% of unit cost.

All hot water tanks presented in Table A.8 are Mitsubishi pre-plumbed cylinders suitable for operation with ASHPs [40]. Table A.9 presents the cost of each tank with respect to a specific heat pump connection, where "N/A" indicates an infeasible connection.

*Table A.8. Parameters for discrete hot water tanks (thermal storage).*

| Label | Model | $\widehat{V}_k$ (m³) | $\hat{\eta}_k^{ch}$ | $\hat{\eta}_k^{disch}$ | $\widehat{T}_k^{min}$ (°C) | Heat Loss (kW) | References |
| --- | --- | --- | --- | --- | --- | --- | --- |
| V | EHPT15X-UKHDW | 0.15 | 0.9 | 0.9 | 45 | 0.048 | [40], [55] |
| S | EHPT17X-UKHDW | 0.17 | 0.9 | 0.9 | 45 | 0.051 | [40], [55] |
| M | EHPT21X-UKHDW | 0.21 | 0.9 | 0.9 | 45 | 0.064 | [40], [55] |
| L | EHPT30X-UKHDW | 0.3 | 0.9 | 0.9 | 45 | 0.087 | [40], [55] |

*Table A.9. Cost of each tank (£) with respect to the ASHP it is connected to. "N/A" indicates that a connection is not feasible. All prices obtained from [55].*

| $\widehat{CU}_{p,k}$ (£) | V | S | M | L |
| --- | --- | --- | --- | --- |
| S1 | 3,924 | 3,979 | 4,033 | N/A |
| M1 | 4,402 | 4,457 | 4,511 | N/A |
| M2 | N/A | N/A | 6,319 | 6,595 |
| L1 | N/A | N/A | 6,319 | 6,595 |

The capacity and COP of heat pumps can be described as a function of water supply temperature and ambient temperature [41]. In this study, the water supply temperature $\widehat{T}^{ws}$ of the ASHP is assumed to be fixed at 55 °C for all selected ASHP units. Note that the periodic raising of the tank water temperature to 60 °C, which is done to combat Legionella growth [56], is not modelled.

For COP, the data points provided in manufacturer datasheets for each ASHP at $\widehat{T}^{ws} = 55$ °C were fitted to a Sigmoid function, as shown below, which eliminated the need for piecewise linearisations.

$$\widehat{COP}_{t,p} = \frac{L_p}{1 + \exp\left(-k_p \cdot \left(\widehat{T}_t^{air} - x_p^0\right)\right)} + b_p \qquad (0A.35)$$

where $L_p$, $k_p$, $x_p^0$, and $b_p$ are fitting parameters. The parameters derived for each ASHP, as well as the respective coefficients of determination ($R^2$) for each curve, are described in Table A.10.

Similarly, data points for capacity for each ASHP at $\widehat{T}^{ws} = 55\,°C$ were fitted to a cubic polynomial, as described below:

$$\widehat{H}_{t,p}^{hp,max} = a_p^{cap}\left(\widehat{T}_t^{air}\right)^2 + b_p^{cap}\left(\widehat{T}_t^{air}\right)^2 + c_p^{cap}\left(\widehat{T}_t^{air}\right) + d_p^{cap} \qquad (0A.36)$$

where $a_p^{cap}$, $b_p^{cap}$, $c_p^{cap}$, and $d_p^{cap}$ are fitting parameters, which are presented in Table A.11 for each ASHP with the respective $R^2$ for each curve.

*Table A.10. Parameters for obtaining COP using a Sigmoid function.*

| Label | Model | $L_p$ | $x_p^0$ | $k_p$ | $b_p$ | $R^2$ |
|---|---|---|---|---|---|---|
| S1 | PUZ-WM50VHA(-BS) | 1.438 | 5.358 | 0.8645 | 1.915 | 0.9874 |
| M1 | PUZ-WM60VAA(-BS) | 445.4 | -1025 | 0.00047 | -272.9 | 0.9883 |
| M2 | PUZ-HWM140VHA(-BS) | 2142 | 300.8 | 0.02195 | -0.6586 | 0.9794 |
| L1 | 2 x PUZ-HWM140VHA(-BS) | 2142 | 300.8 | 0.02195 | -0.6586 | 0.9794 |

*Table A.11. Parameters for obtaining ASHP capacity using a cubic polynomial function.*

| Label | Model | $a_p^{cap}$ | $b_p^{cap}$ | $c_p^{cap}$ | $d_p^{cap}$ | $R^2$ |
|---|---|---|---|---|---|---|
| S1 | PUZ-WM50VHA(-BS) | 0.000409 | -0.01234 | 0.06176 | 5.5 | 0.9001 |
| M1 | PUZ-WM60VAA(-BS) | 0.000966 | -0.0156 | 0.02112 | 7.49 | 0.966 |
| M2 | PUZ-HWM140VHA(-BS) | 0.000598 | -0.0014 | -0.014 | 14.37 | 0.9907 |
| L1 | 2 x PUZ-HWM140VHA(-BS) | 0.001195 | -0.0028 | -0.028 | 28.75 | 0.9907 |

Table A.12 includes parameters for a selection of discrete boilers available from commercial manufacturers.

*Table A.12. Parameters for discrete gas boiler selection.*

| Model | $\widehat{H}_b^{max}$ (kW) | $\hat{\eta}_b$ | $\widehat{CU}_b$ (£) | $\widehat{CI}_b$ (£) | References |
|---|---|---|---|---|---|
| Ideal Logic Combi 24 | 24 | 0.94 | 792 | 1,500 | [57], [58] |
| Potterton Combi Assure 25 | 25 | 0.93 | 797 | 1,563 | [57], [58] |
| Viessmann Combi Vitodens 050-W | 29 | 0.976 | 742 | 1,813 | [57], [58] |
| Viessmann Combi Vitodens 200-W | 32 | 0.98 | 1702 | 2,000 | [57], [58] |

Discrete battery options and respective parameters are presented in Table A.13. Battery installation costs are assumed to be 15% of the unit costs. Operational cost for each battery is assumed to be unit cost divided by lifetime, to account for a full battery replacement within the 20 years, taking battery lifetime as 10 years.

*Table A.13. Parameters for discrete batteries.*

| Label | Model | Cap (kWh) | Max. DOD | Max. SOC | $\widehat{CU}_c$ (£) | $\widehat{CO}_c$ (£/y) | $\hat{\eta}_c^{ch}$ | $\hat{\eta}_c^{disch}$ | $\widehat{CI}_c$ (£) | $E_c^{max}$ (kW) | Ref. |
|---|---|---|---|---|---|---|---|---|---|---|---|
| RESU6.5 | LG Chem RESU6.5 | 6.5 | 0.9 | 1 | 3200 | 160 | 0.97 | 0.97 | 480 | 4.2 | [59], [60] |
| RESU3.3 | LG Chem RESU3.3 | 3.3 | 0.87 | 1 | 2200 | 110 | 0.97 | 0.97 | 330 | 3 | [59], [60] |
| TP2 | Tesla Powerwall 2 | 14 | 0.95 | 1 | 6000 | 300 | 0.95 | 0.95 | 2000 | 5 | [61], [62] |

# Appendix B – Scalars

*Table B.1. Scalar values used in the case studies.*

| Parameter | Value | Unit | Reference |
|---|---|---|---|
| DES lifetime | 20 | years | - |
| Interest rate | 7.5 | % | - |
| Gas purchasing price | 0.02514 | £/kWh | Private correspondence |
| Economy 7 day tariff | 0.18 | £/kWh | [63] |
| Economy 7 night tariff | 0.08 | £/kWh | [63] |
| Export tariff | 0.0503 | £/kWh | [64] |
| Generation tariff (Case 1a & 1b) | 0.1 | £/kWh | Assumed |
| Generation tariff (Case 1c) | 2 | £/kWh | Assumed |
| PV Investment cost | 450 | £/panel | Average from [65] |
| PV Efficiency | 0.18 | - | [66] |
| PV Fixed operational cost | 12.5 | £/kW-yr | [66] |
| Panel area | 1.75 | $m^2$ | Approximation from [65] |
| Panel capacity | 0.25 | kW | Average from [65] |